\documentclass[a4paper, 12pt]{article} 
\usepackage{amssymb,amsmath,latexsym,color,graphicx} 

\newtheorem{dref}{Definition}[section] 
\newtheorem{theo}[dref]{Theorem} 
 \newtheorem{ex}[dref]{Example}

\newenvironment{proofof}{\par\noindent{{\bf Proof} of }}{\hfill$\Box$
\medskip} 
 
\newcommand{\ekv}[2]{\begin{equation}\label{#1}#2\end{equation}}

\newcommand{\no}[1]{(\ref{#1})} 

\title{PT symmetry and Weyl asymptotics} \author{Johannes Sj{\"o}strand \footnote{Ce travail a
    b\'en\'efici\'e d'une aide de l'Agence Nationale de la Recherche
    portant la r\'ef\'erence ANR-08-BLAN-0228-01 ainsi que d'une
    bourse FABER du conseil r\'egional de Bourgogne}\\\small Institut de
  Math\'ematiques de Bourgogne,
  Universit\'e de Bourgogne\\
  \small 9 avenue Alain Savary - BP 47870\\
  \small 21078 Dijon cedex\\ \footnotesize
  johannes.sjostrand@u-bourgogne.fr\\
  \footnotesize and UMR 5584 du CNRS} \date{}
\begin{document}
\maketitle
\centerline{\bf In memory of Leon Ehrenpreis}

\bigskip
\abstract{For a class of PT-symmetric operators with small
  random perturbations, the eigenvalues obey Weyl asymptotics
  with probability close to 1. Consequently, when the principal symbol
is non-real, there are many non-real eigenvalues.}
\tableofcontents

\section{Introduction}\label{in}
\setcounter{equation}{0}

PT-symmetry has been proposed as an alternative for self-adjointness
in quantum physics \cite{BeBoMe99, BeMa10}. Thus for instance, if we
consider a Schr\"odinger operator on ${\bf R}^n$,
\ekv{in.1}{P=-h^2\Delta +V(x),}
the usual assumption of self-adjointness (implying that the potential
$V$ is real valued) can be replaced by that of PT-symmetry:
\ekv{in.2}{V\circ \upsilon =\overline{V},}
where $\upsilon :\, {\bf R}^n\to {\bf R}^n$ is an isometry with $\upsilon
^2= 1\ne\upsilon $. If we introduce the {\bf p}arity operator $U_\upsilon
u(x)=u(\upsilon (x))$ and the {\bf t}ime reversal operator $\Gamma
u=\overline{u}$, then this can be written
\ekv{in.3}
{
[P,U_\upsilon \Gamma ]=0.
} 

Under mild additional technical assumptions it is easy to see that the
spectrum of a PT-symmetric operator is invariant under reflexion in the
real axis. However, in order to build PT-symmetric quantum physics it
seems important that the spectrum be real, so a natural mathematical
question is then to determine when so is the case. Results on reality
and non-reality of the spectrum of PT-symmetric operators can be found
in \cite{Sh02, CaGrSj05, CaGrSj07, BeMa10}.

The purpose of this note is to show that in a probabilistic sense
``most'' non-self-adjoint PT-symmetric operators that are symmetric in the
sense of (\ref{scl.4}), have their
eigenvalues distributed according to the Weyl law and hence many of
their eigenvalues are non-real. As a matter of fact, this will be a
rather easy adaptation of general results on the Weyl asymptotics for
non-self-adjoint operators with small random perturbations
\cite{Ha06a, Ha06b, HaSj08, Bo11, Sj08a, Sj08b, BoSj10}, where
the last three references are the onces that we shall use directly. For
technical reasons we will state our results for elliptic operators on
compact manifolds but it would be easy to adapt the results of
\cite{Sj08a} in order to treat Schr\"odinger operators on ${\bf R}^n$.

The addition of small random perturbations has the effect of
destroying (uniform) analyticity (if the
unperturbated operator has analytic coefficients). A very interesting
question is to give criteria for PT symmetric operators with analytic
coefficients to have real spectrum. 

The plan of the paper is the following: In Section \ref{scl} we treat
the semi-classical case and in Section \ref{la} we treat the case of
large eigenvalues.

\section{The semi-classical case
}\label{scl}
\setcounter{equation}{0}

Let $X$ be a compact smooth manifold of dimension $n$.  Let $\upsilon
:X\to X$ be a smooth involution; $\upsilon ^2=\mathrm{id}$, with $\upsilon
\ne \mathrm{id}$. Fix a smooth positive density $dx$ on $X$ which is
invariant under $\upsilon $ and let us take $L^2$ norms with respect to
$dx$. Let $P$ be a a differential operator on $X$ of order $m\ge 2$
with smooth coefficients so that in local coordinates,
\ekv{scl.1}{ P=\sum_{|\alpha |\le m}a_\alpha
  (x;h)(hD_x)^\alpha, \quad a_\alpha (\cdot ;h)\in C^\infty .}
Here $0<h\ll 1$ is the semi-classical parameter and we assume that 
\ekv{scl.2}
{
a_\alpha (x;h)-a_\alpha (x;0)={\cal O}(h)
}
locally uniformly and similarly for all its derivatives. We also
assume for simplicity that $a_\alpha (x;h)=a_\alpha (x)$ is
independent of $h$ when $|\alpha |=m$.
 Let
$$
p(x,\xi )=\sum_{|\alpha |\le m}a_\alpha (x;0)\xi ^\alpha ,\ \ p_m(x,\xi )=\sum_{|\alpha |= m}a_\alpha (x)\xi ^\alpha 
$$
We assume that $p_m(x,\xi )\ne 0$ on $T^*X\setminus 0$, so that $P$ is
elliptic in the classical sense. We also assume that
\ekv{scl.3}{p_m(T^*X)\ne {\bf C}.} 

\par Assume that $P$ is symmetric, 
\ekv{scl.4}{ P=\Gamma P^*\Gamma =:P^\mathrm{t}.} and that
\ekv{scl.5}{PU=UP^*,\hbox{ where } Uu(x)=U_\upsilon u(x):=u(\upsilon
  (x)),\ \Gamma u(x)=\overline{u(x)}.} This means that $P$ is PT
symmetric: \ekv{scl.6}{[U\Gamma ,P]=0.}  In addition to the
PT-symmetry property (\ref{scl.6}), we have assumed in (\ref{scl.4})
that $P$ is symmetric.

\begin{ex}
$P=-h^2\Delta +V(x)$ on ${\bf T}^n$ where $\Re V$ is even and $\Im V$
is odd, $V(-x)=\overline{V}(x)$. Then $P$ is symmetric and
PT-symmetric with $\upsilon (x)=-x$.
\end{ex}

Let $\widetilde{R}$ be an auxiliary $h$-independent positive elliptic
second order differential operator on $X$ which commutes with
$U$. We also assume that $\widetilde{R}$ is real, or equivalently that
\ekv{scl.6.1}{[\Gamma ,\widetilde{R}]=0.}
 Then $\widetilde{R}$ has an orthonormal basis of real eigenfunctions
$e_j$ such that $Ue_j=(-1)^{k(j)}e_j$ where $k(j)=1$ or
$k(j)=-1$. We say that $e_j$ is even in the first case and odd in
second case. Put $\epsilon _j=e_j$ when $e_j$ is even and
$\epsilon _j=ie_j$ when $e_j$ is odd. Then $\{ \epsilon _j\}$ is also
an orthonormal basis and a linear combination $V=\sum \alpha _j
\epsilon _j$ is PT symmetric iff the coefficients $\alpha _j$ are real:
$U(V)=\overline{V}$. 

In order to formulate our result, we shall follow \cite{Sj08b},
where we treated a situation without any extra symmetry. 

\par Let $\Omega \Subset {\bf C}$ be open, simply connected, not
  entirely contained in $\Sigma (p):=p(T^*X)$. Let
$V_z(t):=\mathrm{vol\,}(\{ \rho \in {\bf R}^{2n};\,|p(\rho )-z|^2 \le
t\} )$. For $\kappa \in ]0,1]$, $z\in \Omega $, we consider the
property that \ekv{int.6.2}{V_z(t)={\cal O}(t^ \kappa ),\ 0\le t \ll
  1.}  Since $r\mapsto p(x,r\xi )$ is a polynomial of degree $m$ in
$r$ with non-vanishing leading coefficient, we see that
(\ref{int.6.2}) holds with $\kappa =1/(2m)$.

By $B_{{\bf R}^d}(0,r)$ we denote the open ball in ${\bf R}^d$ with
center $0$ and radius $r$.
Let $q_\omega $ be a random potential of the form,
\ekv{int.6.3}
{q_\omega (x)=\sum_{0<\mu _k\le L}\alpha _k(\omega )\epsilon
_k(x),\ \alpha (\omega )=(\alpha _k(\omega ))_{0<\mu _k\le L}\in B_{{\bf R}^D}(0, R),}
where $\mu _k>0$ are the square roots of the eigenvalues of
$h^2\widetilde{R}$ so that $h^2\widetilde{R}\epsilon _k=\mu
_k^2 \epsilon _k$. 
We choose $L=L(h)$, $R=R(h)$ in the interval
\ekv{int.6.4}
{\begin{split}h^{\frac{\kappa -3n}{s-\frac{n}{2}-\epsilon }}\ll L\le
    Ch^{-M},&\ \
M\ge \frac{3n-\kappa }{s-\frac{n}{2}-\epsilon },\\\frac{1}{C}h^{-(\frac{n}{2}+\epsilon )M+\kappa -\frac{3n}{2}}\le R\le
 C h^{-\widetilde{M}},&\ \ \widetilde{M}\ge \frac{3n}{2}-\kappa
  +(\frac{n}{2}+\epsilon )M,\end{split}}
for some $\epsilon \in ]0,s-\frac{n}{2}[$, $s>\frac{n}{2}$,
so by Weyl's law for the large eigenvalues of elliptic
self-adjoint operators, the dimension $D$ {in (\ref{int.6.3})} is of the order of magnitude
$(L/h)^n$. We introduce the  small parameter 
$\delta =\tau _0 h^{N_1+n}$, $0<\tau _0\le \sqrt{h}$, where 
\ekv{int.6.4.3}
{
N_1:=\widetilde{M}+sM+\frac{n}{2}.
} 
The randomly perturbed PT symmetric operator is
\ekv{int.6.4.5}
{
P_\delta =P+\delta h^{N_1}q_\omega =:P+\delta Q_\omega .
} 
Here (cf \cite{Sj08a}) the exponent $N_1$ has been chosen so that we have
  uniformly for $h\ll 1$ and $q_\omega $ as above:
$$
\Vert h^{N_1}q_\omega \Vert_{L^\infty }\le {\cal O}(1)h^{-n/2}\Vert
h^{N_1}q_\omega \Vert_{H_h^s}\le {\cal O}(1),
$$
where $H_h^s$ is the natural semi-classical Sobolev space discussed in
Section 2 of $\cite{Sj08b}$ with a norm equivalent to
the standard norm in $H^s$ for each fixed $h>0$.  

\par The random variables $\alpha _j(\omega )$ will have a
joint probability distribution \ekv{int.6.5}{P(d\alpha )=C(h)e^{\Phi
(\alpha ;h)}L(d\alpha ),} where for some $N_4>0$,
\ekv{int.6.6}{ |\nabla _\alpha \Phi |={\cal
O}(h^{-N_4}),} and $L(d\alpha )$ is the
Lebesgue measure. ($C(h)$ is the normalizing constant, 
assuring that the probability of
$B_{{\bf R}^D}(0,R)$ is equal to 1.) 

\par We also need the parameter 
\ekv{int.6.7.5}{\epsilon _0(h)=(h^{\kappa }+h^n\ln 
\frac{1}{h})(\ln \frac{1}{\tau _0}+(\ln \frac{1}{h})^2)} and assume
that $\tau _0=\tau _0(h)$ is not too small, so that $\epsilon _0(h)$ is
small. Recall that $\Omega \Subset {\bf C}$ is open, simply connected, not
entirely contained in $\Sigma (p)$. The main result of this section is:
\begin{theo}\label{int1} Under the assumptions above, let 
$\Gamma \Subset \Omega $ have smooth boundary, let $\kappa \in
]0,1]$ be the parameter in \no{int.6.3}, \no{int.6.4}, \no{int.6.7.5} and assume that 
\no{int.6.2} holds uniformly for $z$ in a
neighborhood of $\partial \Gamma $.  Then there
exists a constant $C>0$ such that for
$C^{-1}\ge r>0$,
$\widetilde{\epsilon }\ge C \epsilon _0(h)$ 
we have with probability 
\ekv{int.6.8}{
\ge 1-\frac{C\epsilon _0(h)}
{rh^{n+\max (n(M+1), N_4+\widetilde{M})}}
e^{-\frac{\widetilde{\epsilon }}{C\epsilon _0(h)}} }
that:
\ekv{int.7}
{\begin{split}
&|
\#(\sigma (P_\delta )\cap \Gamma )-\frac{1}{(2\pi h)^n
}\mathrm{vol\,}(p^{-1}(\Gamma ))
|\le
\\
&
\frac{C}{h^n}\left( \frac{\widetilde{\epsilon }}{r}
+C(r+\ln (\frac{1}{r})\mathrm{vol\,}(p^{-1}(\partial
\Gamma +D(0,r))))
 \right).\end{split}}
Here $\#(\sigma (P_\delta )\cap \Gamma )$ denotes the number of
eigenvalues of $P_\delta $ in $\Gamma $, counted with their algebraic multiplicity.
\end{theo}

\par In the introduction of \cite{Sj08a} there is a discussion about
the choice of parameters which applies here also: Very roughly, if $\tau _0$ is equivalent to
some high power of $h$, then up to some power of $\ln (1/h)$,
$\epsilon _0$ is of the order of magnitude $h^\kappa $. Now choose
$\widetilde{\epsilon }=h^{\kappa -\kappa _0}$ for some $\kappa _0\in
]0,\kappa [$. When $\kappa >1/2$, then the volume in (\ref{int.7}) is
${\cal O}(r^\beta )$ with $\beta =2\kappa -1>0$ and more generally we
may assume that it is ${\cal O}(r^\beta )$ for some $\beta >0$. Then
we choose $r$ to be a suitable power of $h$ and obtain that the right hand side
in (\ref{int.7}) is ${\cal O}(h^{\gamma -n}$) for some $\gamma >0$. With
these choices of the parameters we also see that the probability in
(\ref{int.6.8}) is very close to $1$.

\begin{proofof} Theorem \ref{int1}. We just have to make some small
  modifications in the proof of the main result in \cite{Sj08b} (which
  in turn is a modification of the proof in \cite{Sj08a}) and only
  mention the points where a difference appears. The proof in the
  two cited papers (see also the lecture notes \cite{Sj09}) uses three
  ingredients:
\begin{itemize}
\item[1)] The construction of a special perturbation of the form
  $\delta q_\omega $ with $q_\omega $ as in (\ref{int.6.3}) but with
  $\alpha $ in the {\it complex} ball $B_{{\bf C}^D}(0,R)$ for which
  we have nice lower bounds on the small singular values of $P_\delta $ in
  (\ref{int.6.4.5}), see Proposition 7.3 in \cite{Sj08a}, Proposition
  5.1 in \cite{Sj08b}.
\item[2)] A complex variable argument in the $\alpha $ variables using
  the existence of the special perturbation in step 1), which permits
  to conclude that we have nice lower bounds on a relative determinant
  for $P_\delta -z$, with probability close to 1.
\item[3)] Application of a proposition about the number of zeros of holomorphic
  functions with exponential growth. (See also \cite{Sj10} for an
  improved version of this proposition, not yet fully exploited.)
\end{itemize}

In the present situation we want our special perturbation $\delta
q_\omega (x)$ to be PT-symmetric, that is we want the coefficients
$\alpha $ in (\ref{int.6.3}) to be real. All the parts of the proofs
in step 1 immediately carry over to the case of real $\alpha $ except
the following result which is the basic ingredient in the iterative
process leading to the propositions mentioned above:

\par Let $e_1,...,e_N$ be an ON family in $L^2(X)$ such that 
$$
\Vert \sum_1^N \lambda _je_j\Vert_{H_h^s}\le {\cal O}(1)\Vert \lambda \Vert_{{\bf C}^N}
$$
where the constant ${\cal O}(1)$ is independent of the family and
especially of $N$. Then there exists 
\ekv{scl.7}
{
q=\sum _{0<\mu _j\le L}\alpha _j\epsilon _j,\ \alpha _j\in {\bf C},
}
with $\|\alpha\|_{{\bf C}^D}\le R $ with the parameters as in
(\ref{int.6.4}), such that 
$$
\Vert q\Vert_{H^s_h}\le {\cal
  O}(1)h^{-\frac{n}{2}}NL^{s+\frac{n}{2}+\epsilon }
$$
 and such that the matrix 
$$M_q=(\int q(x)e_j(x)e_k(x))dx)_{1\le j,k\le N}$$ and its singular
values 
$$
\Vert M_q\Vert =s_1(M_q)\ge ...\ge s_N(M_q)
$$ satisfy
$$
\Vert M_q\Vert\le {\cal O}(1)Nh^{-n},
$$
\ekv{scl.8}{
s_k(M_q)\ge h^n/{\cal O}(1),\hbox{ for }1\le k\le N/2.}
(See (6.23), (7.20), (7.23) in \cite{Sj08a}.)

\par Write $q=q_1+iq_2$ where $q_1=\sum (\Re \alpha _j) \epsilon _j$,
$q_2=\sum (\Im \alpha _j) \epsilon _j$, so that $q_1$ and $q_2$ are
PT-symmetric. The upper bounds on $\Vert q\Vert_{H^s_h}$ and on
$\Vert M_q\Vert$ follow from the bound $\Vert \alpha \Vert\le R$
and therefore carry over to $q_j$. Since $M_q=M_{q_1}+iM_{q_2}$ we can
apply the Ky Fan inequalities (\cite{GoKr69}) and get
$$
\frac{h^n}{{\cal O}(1)}\le s_{2k-1}(M_q)\le
s_k(M_{q_1})+s_k(M_{q_{2}}),\ 1\le k\le \frac{N}{4}.
$$
Since the singular values are enumerated in decreasing order, it follows
that for $j$ equal to 1 or 2, we have
\ekv{scl.9}
{
s_k(M_{q_j})\ge \frac{h^n}{2{\cal O}(1)},\ 1\le k\le \frac{N}{4}.
}
this means that step 1 can be carried out and we get a PT symmetric
operator $P_\delta $ as in Proposition
  5.1 in \cite{Sj08b}, the only slight difference is that rather than
  taking $\theta$ in $]0,1/4[$ we have to confine this parameter to
  the smaller interval $]0,1/8[$.

\par Step 2 now follows follows from Remark 8.3 in \cite{Sj08a}, where
the main point is the reality of the coefficients $\alpha _j$ while the
assumption of reality of the basis elements is not necessary, and was
made there only because we had in mind a real perturbation.

\par Step 3 can be carried out without any modifications. 
\end{proofof} 

\section{Weyl asymptotics for large eigenvalues}\label{la}
\setcounter{equation}{0}

Let $P^0$ be an
elliptic differential operator on $X$ of order $m\ge 2$ with smooth
coefficients and with
principal symbol $p_m(x,\xi )$. In local coordinates we get, using
standard multi-index notation,
\ekv{la.1}
{
P^0=\sum_{|\alpha |\le m}a_\alpha ^0(x)D^\alpha ,\quad 
p_m(x,\xi )=\sum_{|\alpha |= m}a_\alpha ^0(x)\xi ^\alpha.
}
Recall that the ellipticity of $P^0$ means that $p_m(x,\xi )\ne 0$ for
$\xi \ne 0$. We assume that
\ekv{la.2}
{
p_m(T^*X)\ne {\bf C}.
}
As before we assume symmetry,
\ekv{la.3}
{
(P^0)^*=\Gamma P^0\Gamma ,
}
and that
\ekv{la.4}
{
P^0U=U(P^0)^*,
}
with $U=U_\upsilon $ as in Section \ref{scl}.

\par Let $\widetilde{R}$ be a reference operator as in and around
(\ref{scl.6.1}) and define $\epsilon _j$ as there. Write
 \ekv{la.5} {
  \widetilde{R}\epsilon _j=(\mu _j^0)^2\epsilon _j,\quad 0<\mu
  _0^0<\mu _1^0\le \mu _2^0\le ...  }
so that $\mu _k=h\mu _k^0$ where $\mu _k$ are given after (\ref{int.6.3}).
 Our randomly perturbed operator
is \ekv{la.6} { P_\omega ^0=P^0+q_\omega ^0(x), } where $\omega $ is the
random parameter and \ekv{la.7} { q_\omega ^0(x)=\sum_{0}^\infty
  \alpha _j^0(\omega )\epsilon _j.  } Here we assume that $\alpha
_j^0(\omega )$ are independent real Gaussian random variables of
variance $\sigma _j^2$ and mean value 0: \ekv{la.8} { \alpha _j^0\sim
  {\cal N}(0,\sigma _j^2), } where \ekv{la.8.5} { (\mu _j^0)^{-\rho
  }e^{-(\mu _j^0)^{\frac{\beta }{M+1}}}\lesssim \sigma _j\lesssim (\mu
  _j^0)^{-\rho }, } \ekv{la.9}
{M=\frac{3n-\frac{1}{2}}{s-\frac{n}{2}-\epsilon },\ 0\le \beta
  <\frac{1}{2},\ \rho >n, } where $s$, $\rho $, $\epsilon $ are fixed
constants such that
$$
\frac{n}{2}<s<\rho -\frac{n}{2},\ 0<\epsilon <s-\frac{n}{2}.
$$

\par Let $H^s(X)$ be the standard Sobolev space of order $s$. As we
saw in \cite{BoSj10} (where the random variables $\alpha _j^0$ were
complex valued), $q_\omega^0 \in H^s(X)$ almost
surely since $s<\rho -\frac{n}{2}$. Hence $q_\omega^0\in L^\infty $
almost surely, implying that $P_\omega ^0$ has purely discrete
spectrum.

\par Consider the function $F(w )=\mathrm{arg\,}p_m(w )$ on
$S^*X$. For given $\theta _0\in S^1\simeq {\bf R}/(2\pi {\bf Z})$,
$N_0\in \dot{{\bf N}}:={\bf N}\setminus \{ 0\}$, we introduce the property $P(\theta _0,N_0)$:
\ekv{la.10}
{
\sum_1^{N_0}|\nabla ^kF(w )|\ne 0\hbox{ on }\{ w \in S^*X;\,
F(w )=\theta _0\}.
}
Notice that if $P(\theta _0,N_0)$ holds, then $P(\theta ,N_0)$ holds
for all $\theta $ in some neighborhood of $\theta _0$. Also notice
that if $X$ is connected and $X$, $p$ are analytic and the analytic
function $F$ is non constant, then $\exists N_0\in \dot{{\bf N}}$ such
that $P(\theta _0,N_0)$ holds for all $\theta _0$.

\par We can now state the main result of this section, which is an
adaptation of the main result of \cite{BoSj10}.

\begin{theo}\label{la1}
Assume that $m\ge 2$. Let $0\le \theta _1\le \theta _2\le 2\pi $ and
assume that $P(\theta _1,N_0)$ and $P(\theta _2,N_0)$ hold for some
$N_0\in\dot{{\bf N}}$. Let $g\in C^\infty ([\theta _1,\theta
_2];]0,\infty [)$ and put 
$$
\Gamma ^g_{\theta _1,\theta _2;0,\lambda }=\{ re^{i\theta } ; \theta
  _1\le \theta \le \theta _2,\ 0\le r\le \lambda g(\theta )\}.
$$
Then for every $\delta \in ]0,\frac{1}{2}-\beta [$ there exists $C>0$ such
that almost surely: $\exists C(\omega )<\infty $ such that for all
$\lambda \in [1,\infty [$:
\ekv{la.11}
{\begin{split}
&|\#(\sigma (P_\omega ^0)\cap \Gamma _{\theta _1,\theta _2;0,\lambda }^g)
-\frac{1}{(2\pi )^n}\mathrm{vol\,}p_m^{-1}(\Gamma ^g_{\theta _1,\theta
  _2;0,\lambda })
|
\\
&\le C(\omega )+C\lambda ^{\frac{n}{m}-\frac{1}{m}(\frac{1}{2}-\beta -\delta
    )\frac{1}{N_0+1}}.\end{split}}
\end{theo}

The proof actually allows to have almost surely a simultaneous
conclusion for a whole family of $\theta _1,\theta _2,g$:

\begin{theo}\label{la2}
Assume that $m\ge 2$. Let $\Theta $ be a compact subset of $[0,2\pi ]$. Let
$N_0\in {\bf N}$ and assume that $P(\theta ,N_0)$ holds uniformly for
$\theta \in \Theta $. Let ${\cal G}$ be a subset of $\{(g,\theta
_1,\theta _2);\ \theta _j\in \Theta, \theta _1\le \theta _2,\ g\in 
C^\infty ([\theta _1,\theta
_2];]0,\infty [)
\}$ with the property that $g$ and $1/g$ are uniformly bounded in $C^\infty ([\theta _1,\theta
_2];]0,\infty [)$ when $(g,\theta _1,\theta _2)$ varies in ${\cal
  G}$. Then for every $\delta \in ]0,\frac{1}{2}-\beta [$ 
there exists $C>0$ such
that almost surely: $\exists C(\omega )<\infty $ such that for all
$\lambda \in [1,\infty [$ and all $(g,\theta _1,\theta _2)\in {\cal
  G}$, we have the estimate (\ref{la.11}).
\end{theo}

The condition (\ref{la.8.5}) allows us to choose $\sigma _j$ decaying
faster than any negative power of $\mu _j^0$. Then from the discussion
below, it will follow that $q_\omega (x)$ is almost surely a smooth
function. A rough and somewhat intuitive interpretation of Theorem
\ref{la2} is then that for almost every PT symmetric elliptic operator
of order $\ge 2$ with smooth coefficients on a compact manifold which
satisfies the conditions (\ref{la.2}), (\ref{la.3}), (\ref{la.4}), the
large eigenvalues distribute according to Weyl's law in sectors with
limiting directions that satisfy a weak non-degeneracy condition.

\medskip
\begin{proofof} Theorem \ref{la1}. 
As already mentioned, the theorem is a variant of Theorem 1.1 in
\cite{BoSj10}. The difference is just that we now use real
random variables in the perturbation $q_\omega ^0$ in order to assure
the PT-symmetry while in \cite{BoSj10} they were complex. The proof in
\cite{BoSj10} used a reduction to the semi-classical case where the
main result of \cite{Sj08b} could be applied. 
The proof of Theorem \ref{la1} is an immediate modification of that
proof, where we replace the main result in \cite{Sj08b} by Theorem
\ref{int1}. The only point where the use of real Gaussian random
variables in stead of complex ones causes a slight change is the use of
(4.10) in \cite{BoSj10} that was established in \cite{Bo08}, where we
have to replace the denominator 2 by 4 in the case of real random
variables. That was also proved by Bordeaux Montrieux in \cite{Bo08},
Proposition 2.5.4.
\end{proofof}

\end{document}